\documentclass{article}

\usepackage{amsmath,amssymb,amsthm}

\newtheorem{thm}{Theorem}
\newtheorem{cor}{Corollary}
\newtheorem{lem}{Lemma}
\newtheorem{prop}{Proposition}

\theoremstyle{remark}
\newtheorem{rem}{Remark}

\theoremstyle{definition}
\newtheorem{defin}{Definition}

\title{Necessary optimality conditions
for the calculus of variations on time scales\footnote{This work
is part of the first author's PhD project.}}

\author{Rui A. C. Ferreira\\
\texttt{ruiacferreira@yahoo.com}
\and Delfim F. M. Torres\\
\texttt{delfim@mat.ua.pt}}

\date{Department of Mathematics\\
University of Aveiro\\
3810-193 Aveiro, Portugal}

\begin{document}

\maketitle

\begin{abstract}
We study more general variational problems on time scales.
Previous results are generalized by proving necessary optimality
conditions for (i) variational problems involving delta
derivatives of more than the first order, and (ii) problems of the
calculus of variations with delta-differential side conditions
(Lagrange problem of the calculus of variations on time scales).

\bigskip

\noindent \textbf{Keywords:} time scales, $\Delta$-variational
calculus, higher-order $\Delta$-derivatives, higher-order
Euler-Lagrange $\Delta$-equations, Lagrange problem on time
scales, normal and abnormal $\Delta$-extremals.

\bigskip

\noindent \textbf{2000 Mathematics Subject Classification:} 49K05, 39A12.
\end{abstract}


\section{Introduction}

The theory of time scales is a relatively new area, that unify and
generalize difference and differential equations \cite{livro}. It
was initiated by Stefan Hilger in the nineties of the XX century
\cite{Hilger90,Hilger97}, and is now subject of strong current
research in many different fields in which dynamic processes can
be described with discrete or continuous models \cite{Agarwal}.

The calculus of variations on time scales was introduced by Bohner
\cite{CD:Bohner:2004} and by Hilscher and Zeidan \cite{zeidan},
and appears to have many opportunities for application in
economics \cite{Atici06}. In all those works, necessary optimality
conditions are only obtained for the basic (simplest) problem of
the calculus of variations on time scales: in
\cite{Atici06,CD:Bohner:2004} for the basic problem with fixed
endpoints, in \cite{zeidan} for the basic problem with general
(jointly varying) endpoints. Having in mind the classical setting
(situation when the time scale $\mathbb{T}$ is either $\mathbb{R}$
or $\mathbb{Z}$ -- see \textrm{e.g.} \cite{GelfandFomin,Brunt} and
\cite{zeidan2,Logan}, respectively), one suspects that the
Euler-Lagrange equations in \cite{Atici06,CD:Bohner:2004,zeidan}
are easily generalized for problems with higher-order delta
derivatives. This is not exactly the case, even beginning with the
formulation of the problem.

The basic problem of the calculus of variations on time scales is
defined (\textrm{cf.} \cite{CD:Bohner:2004,zeidan}, see
\S\ref{sec:Prel} below for the meaning of the $\Delta$-derivative
and $\Delta$-integral) as
\begin{equation}
\label{eq:EL:B} \mathcal{L}[y(\cdot)]
=\int_{a}^{b}L(t,y^\sigma(t),y^\Delta(t))\Delta t\longrightarrow
\min, \quad (y(a)=y_a)\, , (y(b)=y_b) \, ,
\end{equation}
with $L: \mathbb{T} \times \mathbb{R}^n \times \mathbb{R}^n
\rightarrow \mathbb{R}$, $(y,u) \rightarrow L(t,y,u)$ a
$C^2$-function for each $t$, and where we are using parentheses
around the endpoint conditions as a notation to mean that the
conditions may or may not be present: the case with fixed boundary
conditions $y(a)=y_a$ and $y(b)=y_b$ is studied in
\cite{CD:Bohner:2004}, for admissible functions $y(\cdot)$
belonging to $C^1_{rd}\left(\mathbb{T};\mathbb{R}^n\right)$
($rd$-continuously $\Delta$-differentiable functions); general
boundary conditions of the type $f(y(a),y(b))=0$, which include
the case $y(a)$ or $y(b)$ free, and over admissible functions in
the wider class $C^1_{prd}\left(\mathbb{T};\mathbb{R}^n\right)$
(piecewise $rd$-continuously $\Delta$-differentiable functions),
are considered in \cite{zeidan}. One question immediately comes to
mind. Why is the basic problem on time scales defined as
\eqref{eq:EL:B} and not as
\begin{equation}
\label{eq:EL:BSS} \mathcal{L}[y(\cdot)] =\int_{a}^{b}
L(t,y(t),y^\Delta(t))\Delta t\longrightarrow \min, \quad
(y(a)=y_a)\, , (y(b)=y_b) \, .
\end{equation}
The answer is simple: compared with \eqref{eq:EL:BSS}, definition
\eqref{eq:EL:B} simplifies the Euler-Lagrange equation, in the
sense that makes it similar to the classical context. The reader
is invited to compare the Euler-Lagrange condition
\eqref{eq:EL:Boh} of problem \eqref{eq:EL:B} and the
Euler-Lagrange condition \eqref{minha:E-L} of problem
\eqref{eq:EL:BSS}, with the classical expression (on the time
scale $\mathbb{T} = \mathbb{R}$):
\begin{equation*}
\frac{d}{dt} L_{y'}(t,y_\ast(t),y_\ast'(t))
=L_{y}(t,y_\ast(t),y_\ast'(t)),\ t\in[a,b] \, .
\end{equation*}
It turns out that problems \eqref{eq:EL:B} and \eqref{eq:EL:BSS}
are equivalent: as far as we are assuming $y(\cdot)$ to be
$\Delta$-differentiable, then
$y(t)=y^{\sigma}(t)-\mu(t)y^{\Delta}(t)$ and (i) any problem
\eqref{eq:EL:B} can be written in the form \eqref{eq:EL:BSS}, (ii)
any problem \eqref{eq:EL:BSS} can be written in the form
\eqref{eq:EL:B}. We claim, however, that the formulation
\eqref{eq:EL:BSS} we are promoting here is more natural and
convenient. An advantage of our formulation \eqref{eq:EL:BSS} with
respect to \eqref{eq:EL:B} is that it makes clear how to
generalize the basic problem on time scales to the case of a
Lagrangian $L$ containing delta derivatives of $y(\cdot)$ up to an
order $r$, $r \ge 1$. The higher-order problem will be naturally
defined as
$$\mathcal{L}[y(\cdot)]=\int_{a}^{\rho^{r-1}(b)}
L(t,y(t),y^\Delta(t),\ldots,y^{\Delta^r}(t))\Delta
t\longrightarrow\min,$$
\begin{align}
\left(y(a)=y_a^0\right),& \ \left(y\left(\rho^{r-1}(b)\right)=y_b^0\right), \label{problema } \\
&\vdots\nonumber \\
\left(y^{\Delta^{r-1}}(a)=y_a^{r-1}\right),&\
\left(y^{\Delta^{r-1}}\left(\rho^{r-1}(b)\right)=y_b^{r-1}\right),\nonumber
\end{align}
where $y^{\Delta^i}(t)\in\mathbb{R}^n,\ i\in\{0,\ldots,r\}$,
$y^{\Delta^0}=y$, and $n,\ r\in\mathbb{N}$ (assumptions on the
data of the problem will be specified later, in
Section~\ref{sec:mainResults}). One of the new results in this
paper is a necessary optimality condition in \emph{delta integral
form} for problem \eqref{problema } (Theorem~\ref{thm:HO:E-L:TS}).
It is obtained using the interplay of problems \eqref{eq:EL:B} and
\eqref{eq:EL:BSS} in order to deal with more general optimal
control problems \eqref{eq:PrbCO}.

The paper is organized as follows. In Section~\ref{sec:Prel} we
give a brief introduction to time scales and recall the main
results of the calculus of variations on this general setting. Our
contributions are found in Section~\ref{sec:mainResults}. We start
in \S\ref{subsec:basic} by proving the Euler-Lagrange equation and
transversality conditions (natural boundary conditions -- $y(a)$
or/and $y(b)$ free) for the basic problem \eqref{eq:EL:BSS}
(Theorem~\ref{thm:1}). As a corollary, the Euler-Lagrange equation
in \cite{CD:Bohner:2004} and \cite{zeidan} for \eqref{eq:EL:B} is
obtained. Regarding the natural boundary conditions, the one which
appears when $y(a)$ is free turns out to be simpler and more close
in aspect to the classical condition
$L_{y'}(a,y_\ast(a),y_\ast'(a)) = 0$ for problem \eqref{eq:EL:B}
than to \eqref{eq:EL:BSS}---compare condition \eqref{eq:transv:a}
for problem \eqref{eq:EL:BSS} with the correspondent condition
\eqref{eq:trv:prbBcsig:a} for problem \eqref{eq:EL:B}; but the
inverse situation happens when $y(b)$ is free---compare condition
\eqref{eq:trv:prbBcsig:b} for problem \eqref{eq:EL:B} with the
correspondent condition \eqref{eq:transv:b} for \eqref{eq:EL:BSS},
this last being simpler and more close in aspect to the classical
expression $L_{y'}(b,y_\ast(b),y_\ast'(b)) = 0$ valid on the time
scale $\mathbb{T} = \mathbb{R}$. In \S\ref{subsec:Lag} we
formulate  a more general optimal control problem \eqref{eq:PrbCO}
on time scales, proving respective necessary optimality conditions
in Hamiltonian form (Theorem~\ref{thm:PMP}). As corollaries, we
obtain a Lagrange multiplier rule on time-scales
(Corollary~\ref{cor:LagMultRule}), and in \S\ref{subsec:HO} the
Euler-Lagrange equation for the problem of the calculus of
variations with higher order delta derivatives
(Theorem~\ref{thm:HO:E-L:TS}). Finally, as an illustrative
example, we consider in \S\ref{subsec:appl} a discrete time scale
and obtain the well-known Euler-Lagrange equation in delta
differentiated form.

All the results obtained in this paper can be extended: (i) to
nabla derivatives (see \cite[\S 8.4]{livro}) with the appropriate
modifications and as done in \cite{Atici06} for the simplest
functional; (ii) to more general classes of admissible functions
and to problems with more general boundary conditions, as done in
\cite{zeidan} for the simplest functional of the calculus of
variations on time scales.


\section{Time scales and previous results}
\label{sec:Prel}

We begin by recalling the main definitions and properties of time
scales (\textrm{cf.} \cite{Agarwal,livro,Hilger90,Hilger97} and
references therein).

A nonempty closed subset of $\mathbb{R}$ is called a \emph{Time
Scale} and is denoted by $\mathbb{T}$.

The \emph{forward jump operator}
$\sigma:\mathbb{T}\rightarrow\mathbb{T}$ is defined by
$$\sigma(t)=\inf{\{s\in\mathbb{T}:s>t}\},\mbox{ for all $t\in\mathbb{T}$},$$
while the \emph{backward jump operator}
$\rho:\mathbb{T}\rightarrow\mathbb{T}$ is defined by
$$\rho(t)=\sup{\{s\in\mathbb{T}:s<t}\},\mbox{ for all
$t\in\mathbb{T}$},$$ with $\inf\emptyset=\sup\mathbb{T}$
(\textrm{i.e.}, $\sigma(M)=M$ if $\mathbb{T}$ has a maximum $M$)
and $\sup\emptyset=\inf\mathbb{T}$ (\textrm{i.e.}, $\rho(m)=m$ if
$\mathbb{T}$ has a minimum $m$).

A point $t\in\mathbb{T}$ is called \emph{right-dense},
\emph{right-scattered}, \emph{left-dense} and
\emph{left-scattered} if $\sigma(t)=t$, $\sigma(t)>t$, $\rho(t)=t$
and $\rho(t)<t$, respectively.

Throughout the text we let $\mathbb{T}=[a,b]\cap\mathbb{T}_{0}$
with $a<b$ and $\mathbb{T}_0$ a time scale. We define
$\mathbb{T}^k=\mathbb{T}\backslash(\rho(b),b]$,
$\mathbb{T}^{k^2}=\left(\mathbb{T}^k\right)^k$ and more generally
$\mathbb{T}^{k^n}=\left(\mathbb{T}^{k^{n-1}}\right)^k$, for
$n\in\mathbb{N}$. The following standard notation is used for
$\sigma$ (and $\rho$): $\sigma^0(t) = t$, $\sigma^n(t) = (\sigma
\circ \sigma^{n-1})(t)$, $n \in \mathbb{N}$.

The \emph{graininess function}
$\mu:\mathbb{T}\rightarrow[0,\infty)$ is defined by
$$\mu(t)=\sigma(t)-t,\mbox{ for all $t\in\mathbb{T}$}.$$

We say that a function $f:\mathbb{T}\rightarrow\mathbb{R}$ is
\emph{delta differentiable} at $t\in\mathbb{T}^k$ if there is a
number $f^{\Delta}(t)$ such that for all $\varepsilon>0$ there
exists a neighborhood $U$ of $t$ (\textrm{i.e.},
$U=(t-\delta,t+\delta)\cap\mathbb{T}$ for some $\delta>0$) such
that
$$|f(\sigma(t))-f(s)-f^{\Delta}(t)(\sigma(t)-s)|
\leq\varepsilon|\sigma(t)-s|,\mbox{ for all $s\in U$}.$$
We call $f^{\Delta}(t)$ the \emph{delta derivative} of $f$ at $t$.

Now, we define the $r^{th}-$\emph{delta derivative}
($r\in\mathbb{N}$) of $f$ to be the function
$f^{\Delta^r}:\mathbb{T}^{k^r}\rightarrow\mathbb{R}$, provided
$f^{\Delta^{r-1}}$ is delta differentiable on $\mathbb{T}^{k^r}$.

For delta differentiable $f$ and $g$, the next formulas hold:

\begin{align}
f^\sigma(t)&=f(t)+\mu(t)f^\Delta(t)\label{transfor}\\
(fg)^\Delta(t)&=f^\Delta(t)g^\sigma(t)+f(t)g^\Delta(t)\nonumber\\
&=f^\Delta(t)g(t)+f^\sigma(t)g^\Delta(t)\nonumber,
\end{align}

where we abbreviate $f\circ\sigma$ by $f^\sigma$.

Next, a function $f:\mathbb{T}\rightarrow\mathbb{R}$ is called
\emph{rd-continuous} if it is continuous at right-dense points and
if its left-sided limit exists at left-dense points. We denote the
set of all rd-continuous functions by C$_{\textrm{rd}}$ or
C$_{\textrm{rd}}[\mathbb{T}]$ and the set of all delta
differentiable functions with rd-continuous derivative by
C$_{\textrm{rd}}^1$ or C$_{\textrm{rd}}^1[\mathbb{T}]$.

It is known that rd-continuous functions possess an
\emph{antiderivative}, \textrm{i.e.}, there exists a function $F$
with $F^\Delta=f$, and in this case an \emph{integral} is defined
by $\int_{a}^{b}f(t)\Delta t=F(b)-F(a)$. It satisfies
\begin{equation}
\label{sigma}
\int_t^{\sigma(t)}f(\tau)\Delta\tau=\mu(t)f(t) \, .
\end{equation}

We now present some useful properties of the delta integral:

\begin{lem}
\label{integracao:partes}
If $a,b\in\mathbb{T}$ and $f,g\in$C$_{\textrm{rd}}$, then
\begin{enumerate}

 \item$\int_{a}^{b}f(\sigma(t))g^{\Delta}(t)\Delta t
 =\left[(fg)(t)\right]_{t=a}^{t=b}-\int_{a}^{b}f^{\Delta}(t)g(t)\Delta t$.

\item $\int_{a}^{b}f(t)g^{\Delta}(t)\Delta t
=\left[(fg)(t)\right]_{t=a}^{t=b}-\int_{a}^{b}f^{\Delta}(t)g(\sigma(t))\Delta t$.

\end{enumerate}
\end{lem}

The main result of the calculus of variations on time scales is
given by the following necessary optimality condition for problem
\eqref{eq:EL:B}.

\begin{thm}[\cite{CD:Bohner:2004}]
\label{Th:B:EL-CV} If $y_\ast$ is a weak local minimizer
(\textrm{cf.} \S\ref{sec:mainResults}) of the problem
\begin{gather*}
\mathcal{L}[y(\cdot)]=\int_{a}^{b}L(t,y^\sigma(t),y^\Delta(t))\Delta t \longrightarrow \min\\
y(\cdot) \in C_{\textrm{rd}}^1[\mathbb{T}]\\
y(a)=y_a, \quad y(b)=y_b,
\end{gather*}
then the Euler-Lagrange equation
\begin{equation}
\label{eq:EL:Boh}
L_{y^\Delta}^\Delta(t,y^\sigma_\ast(t),y_\ast^\Delta(t))
=L_{y^\sigma}(t,y^\sigma_\ast(t),y_\ast^\Delta(t)),\
t\in\mathbb{T}^{k^2}
\end{equation}
holds.
\end{thm}

Main ingredients to prove Theorem~\ref{Th:B:EL-CV} are item~1 of
Lemma~\ref{integracao:partes} and the Dubois-Reymond lemma:

\begin{lem}[\cite{CD:Bohner:2004}]
\label{lem:DR} Let $g\in C_{\textrm{rd}}$,
$g:[a,b]^k\rightarrow\mathbb{R}^n$. Then,
$$\int_{a}^{b}g(t) \cdot \eta^\Delta(t)\Delta t=0  \quad
\mbox{for all $\eta\in C_{\textrm{rd}}^1$ with
$\eta(a)=\eta(b)=0$}$$ if and only if $$g(t)=c \mbox{ on $[a,b]^k$
for some $c\in\mathbb{R}^n$}.$$
\end{lem}


\section{Main results}
\label{sec:mainResults}

Assume that the Lagrangian $L(t,u_0(t),u_1(t),\ldots,u_r(t))$
($r\geq 1$) is a $\mathrm{C}^{r+1}$ function of
$(u_0(t),u_1(t),\ldots,u_r(t))$ for each $t\in\mathbb{T}$. Let
$y\in\mathrm{C}_{rd}^r[\mathbb{T}]$, where
$$\mathrm{C}_{rd}^r[\mathbb{T}]
=\left\{y:\mathbb{T}^{k^r}\rightarrow\mathbb{R}^n : y^{\Delta^r}\
\mbox{is $rd$-continuous on}\ \mathbb{T}^{k^r}\right\} \, .$$

We want to minimize the functional $\mathcal{L}$ of problem
\eqref{problema }. For this, we say that
$y_\ast\in\mathrm{C}_{rd}^r[\mathbb{T}]$ is a \emph{weak local
minimizer} for the variational problem \eqref{problema } provided
there exists $\delta>0$ such that
$\mathcal{L}[y_\ast]\leq\mathcal{L}[y]$ for all
$y\in\textrm{C}_{rd}^r[\mathbb{T}]$ satisfying the constraints in
\eqref{problema } and $\|y-y_\ast\|_{r,\infty}<\delta$, where
$$||y||_{r,\infty} := \sum_{i=0}^{r} \left\|y^{\Delta^i}\right\|_{\infty},$$
with $y^{\Delta^0} = y$ and $||y||_{\infty}:= \sup_{t
\in\mathbb{T}^{k^r}} |y(t)|$.


\subsection{The basic problem on time scales}
\label{subsec:basic}

We start by proving the necessary optimality condition for the
simplest variational problem ($r = 1$):

\begin{equation}
\label{P1}
\begin{gathered}
\mathcal{L}[y(\cdot)]=\int_{a}^{b}L(t,y(t),y^\Delta(t))\Delta t \longrightarrow \min \\
y(\cdot) \in C_{\textrm{rd}}^1[\mathbb{T}]\\
\left(y(a)=y_a\right), \quad \left(y(b)=y_b\right) \, .
\end{gathered}
\end{equation}

\begin{rem}
\label{rem:3points} We are assuming in problem \eqref{P1} that the
time scale $\mathbb{T}$ has at least 3 points. Indeed, for the
delta-integral to be defined we need at least 2 points. Assume
that the time scale has only two points: $\mathbb{T} = \{a,b\}$,
with $b=\sigma(a)$. Then,
$\int_{a}^{\sigma(a)}L(t,y(t),y^\Delta(t))\Delta t = \mu(a)
L(a,y(a),y^\Delta(a))$. In the case both $y(a)$ and $y(\sigma(a))$
are fixed, since $y^\Delta(a) = \frac{y(\sigma(a))-y(a)}{\mu(a)}$,
then $\mathcal{L}[y(\cdot)]$ would be a constant for every
admissible function $y(\cdot)$ (there would be nothing to minimize
and problem \eqref{P1} would be trivial). Similarly, for
\eqref{problema } we assume the time scale to have at least $2 r +
1$ points (see Remark~\ref{rem:Pneeds2rp1points}).
\end{rem}

\begin{thm}
\label{thm:1} If $y_\ast$ is a weak local minimizer of \eqref{P1}
(problem \eqref{problema }  with $r=1$), then the Euler-Lagrange
equation in $\Delta$-integral form
\begin{equation}\label{eulerint}
L_{y^\Delta}(t,y_\ast(t),y_\ast^\Delta(t))
=\int_a^{\sigma(t)}L_{y}(\xi,y_\ast(\xi),y_\ast^\Delta(\xi))\Delta\xi+c
\end{equation}
holds $\forall t\in\mathbb{T}^k$ and some $c\in\mathbb{R}^n$.
Moreover, if the initial condition $y(a)=y_a$ is not present
($y(a)$ is free), then the supplementary condition
\begin{equation}
\label{eq:transv:a} L_{y^\Delta}(a,y_\ast(a),y_\ast^\Delta(a))
-\mu(a)L_{y}(a,y_\ast(a),y_\ast^\Delta(a)) = 0
\end{equation}
holds; if $y(b)=y_b$ is not present ($y(b)$ is free), then
\begin{equation}
\label{eq:transv:b}
L_{y^\Delta}(\rho(b),y_\ast(\rho(b)),y_\ast^\Delta(\rho(b)))= 0\,
.
\end{equation}
\end{thm}

\begin{rem}
For the time scale $\mathbb{T} = \mathbb{R}$ equalities
\eqref{eq:transv:a} and \eqref{eq:transv:b} give, respectively,
the well-known \emph{natural boundary conditions}
$L_{y'}(a,y_\ast(a),y_\ast'(a)) = 0$ and
$L_{y'}(b,y_\ast(b),y_\ast'(b))= 0$.
\end{rem}

\begin{proof}
Suppose that $y_\ast(\cdot)$ is a weak local minimizer of
$\mathcal{L}[\cdot]$. Let
$\eta(\cdot)\in$\textrm{C}$_{\textrm{rd}}^1$ and define
$\Phi:\mathbb{R}\rightarrow\mathbb{R}$ by
$$\Phi(\varepsilon)=\mathcal{L}[y_\ast(\cdot)+\varepsilon\eta(\cdot)].$$
This function has a minimum at $\varepsilon=0$, so we must have
$\Phi'(0)=0$. Applying the delta-integral properties and the
integration by parts formula 2 (second item in
Lemma~\ref{integracao:partes}), we have
\begin{equation}
\label{eq:1:5:2}
\begin{aligned}
0&=\Phi'(0)\\
 &=\int_{a}^{b}[L_{y}(t,y_\ast(t),y_\ast^\Delta(t)) \cdot \eta(t)
 +L_{y^\Delta}(t,y_\ast(t),y_\ast^\Delta(t)) \cdot \eta^\Delta(t)]\Delta t\\
 &=\left.\int_a^t L_{y}(t,y_\ast(t),y_\ast^\Delta(t)\Delta t
 \cdot \eta(t)\right|_{t=a}^{t=b}\\
 &\quad -\int_{a}^{b}\left[\int_{a}^{\sigma(t)}L_{y}(\xi,y_\ast(\xi),
 y_\ast^\Delta(\xi))\Delta\xi \cdot \eta^\Delta(t)
 -L_{y^\Delta}(t,y_\ast(t),y_\ast^\Delta(t)) \cdot \eta^\Delta(t)\right]\Delta t \, .
\end{aligned}
\end{equation}
Let us limit the set of all delta-differentiable functions
$\eta(\cdot)$ with $rd$-continuous derivatives to those which
satisfy the condition $\eta(a)=\eta(b)=0$ (this condition is
satisfied by all the admissible variations $\eta(\cdot)$ in the
case both $y(a)=y_a$ and $y(b)=y_b$ are fixed). For these
functions we have
\begin{equation*}
 \int_{a}^{b}\left[L_{y^\Delta}(t,y_\ast(t),y_\ast^\Delta(t))
 -\int_{a}^{\sigma(t)}L_{y}(\xi,y_\ast(\xi),y_\ast^\Delta(\xi))
 \Delta\xi\right] \cdot \eta^\Delta(t)\Delta t = 0 \, .
\end{equation*}
Therefore, by the lemma of Dubois-Reymond (Lemma~\ref{lem:DR}),
there exists a constant $c\in\mathbb{R}^n$ such that \eqref{eulerint} holds:
\begin{equation}
\label{eq:1:5:3} L_{y^\Delta}(t,y_\ast(t),y_\ast^\Delta(t))
-\int_{a}^{\sigma(t)}L_{y}(\xi,y_\ast(\xi),y_\ast^\Delta(\xi))\Delta\xi=c
\, ,
\end{equation}
for all $t\in\mathbb{T}^k$. Because of \eqref{eq:1:5:3},
condition \eqref{eq:1:5:2} simplifies to
\begin{equation*}
\left.\int_a^t L_{y}(t,y_\ast(t),y_\ast^\Delta(t)\Delta t \cdot
\eta(t)\right|_{t=a}^{t=b}+\left.
c\cdot\eta(t)\right|_{t=a}^{t=b}=0,
\end{equation*}
for any admissible $\eta(\cdot)$. If $y(a) = y_a$ is not present
in problem \eqref{P1} (so that $\eta(a)$ need not to be zero),
taking $\eta(t) = t-b$ we find that $c = 0$; if $y(b) = y_b$ is
not present, taking $\eta(t) = t-a$ we find that $\int_a^b
L_{y}(t,y_\ast(t),y_\ast^\Delta(t) = 0$. Applying these two
conditions to \eqref{eq:1:5:3} and having in mind formula
\eqref{sigma}, we may state that
\begin{multline*}
L_{y^\Delta}(a,y_\ast(a),y_\ast^\Delta(a))
-\int_{a}^{\sigma(a)}L_{y}(\xi,y_\ast(\xi),y_\ast^\Delta(\xi))\Delta\xi=0\\
\Leftrightarrow L_{y^\Delta}(a,y_\ast(a),y_\ast^\Delta(a))
-\mu(a)L_{y}(a,y_\ast(a),y_\ast^\Delta(a))=0,
\end{multline*}
and (note that $\sigma(\rho(b))=b$)
\begin{multline*}
L_{y^\Delta}(\rho(b),y_\ast(\rho(b)),y_\ast^\Delta(\rho(b)))
-\int_{a}^{b}L_{y}(\xi,y_\ast(\xi),y_\ast^\Delta(\xi))\Delta\xi=0\\
\Leftrightarrow
L_{y^\Delta}(\rho(b),y_\ast(\rho(b)),y_\ast^\Delta(\rho(b)))=0.
\end{multline*}
\end{proof}

\begin{rem}\label{rem1}
Since $\sigma(t)\geq t,\ \forall t\in\mathbb{T}$, we must have
\begin{multline*}
L_{y^\Delta}(t,y_\ast(t),y_\ast^\Delta(t))
-\int_{a}^{\sigma(t)}L_{y}(\xi,y_\ast(\xi),y_\ast^\Delta(\xi))\Delta\xi=c\\
\Leftrightarrow L_{y^\Delta}(t,y_\ast(t),y_\ast^\Delta(t))
-\mu(t)L_{y}(t,y_\ast(t),y_\ast^\Delta(t)) \\
=\int_{a}^{t}L_{y}(\xi,y_\ast(\xi),y_\ast^\Delta(\xi))\Delta\xi+c,
\end{multline*}
by formula \eqref{sigma}. Delta differentiating both sides, we
obtain
\begin{multline}
\left(L_{y^\Delta}(t,y_\ast(t),y_\ast^\Delta(t))
-\mu(t)L_{y}(t,y_\ast(t),y_\ast^\Delta(t))\right)^\Delta \\
=L_{y}(t,y_\ast(t),y_\ast^\Delta(t)),\
t\in\mathbb{T}^{k^2}.\label{minha:E-L}
\end{multline}
Note that we can't expand the left hand side of this last
equation, because we are not assuming that $\mu(t)$ is delta
differentiable. In fact, generally $\mu(t)$ is not delta
differentiable (see example 1.55, page 21 of \cite{livro}). We say
that \eqref{minha:E-L} is the Euler-Lagrange equation for problem
\eqref{P1} in the \emph{delta differentiated} form.
\end{rem}

As mentioned in the introduction, the formulations of the problems
of the calculus of variations on time scales with
``$\left(t,y^{\sigma}(t),y^\Delta(t)\right)$'' and with
``$\left(t,y(t),y^\Delta(t)\right)$'' are equivalent. It is
trivial to derive previous Euler-Lagrange equation
\eqref{eq:EL:Boh} from our equation \eqref{minha:E-L} and the
other way around (one can derive \eqref{minha:E-L} directly from
\eqref{eq:EL:Boh}).

\begin{cor}
If $y_\ast \in C_{\textrm{rd}}^1[\mathbb{T}]$ is a weak local
minimizer of
$$\mathcal{L}[y(\cdot)]=\int_{a}^{b}L(t,y^\sigma(t),y^\Delta(t))\Delta
t\, , \quad \mbox{$\left(y(a)=y_a\right)$,
$\left(y(b)=y_b\right)$},$$ then the Euler-Lagrange equation
\eqref{eq:EL:Boh} holds. If $y(a)$ is free, then the extra
transversality condition (natural boundary condition)
\begin{equation}
\label{eq:trv:prbBcsig:a}
L_{y^\Delta}(a,y_\ast^\sigma(a),y_\ast^\Delta(a)) = 0
\end{equation}
holds; if $y(b)$ is free, then
\begin{equation}
\label{eq:trv:prbBcsig:b}
L_{y^\sigma}(\rho(b),y_\ast^\sigma(\rho(b)),y_\ast^\Delta(\rho(b)))\mu(\rho(b))
+L_{y^\Delta}(\rho(b),y_\ast^\sigma(\rho(b)),y_\ast^\Delta(\rho(b)))
= 0\, .
\end{equation}
\end{cor}

\begin{proof}
Since $y(t)$ is delta differentiable, then \eqref{transfor} holds.
This permits us to write
$$L(t,y^\sigma(t),y^\Delta(t))=L(t,y(t)+\mu(t)y^\Delta(t),y^\Delta(t))
=F(t,y(t),y^\Delta(t)).$$ Applying equation \eqref{minha:E-L} to
the functional $F$ we obtain
$$\left(F_{y^\Delta}(t,y(t),y^\Delta(t))
-\mu(t)F_{y}(t,y(t),y^\Delta(t))\right)^\Delta=F_{y}(t,y(t),y^\Delta(t)).$$
But
\begin{align}
F_{y}(t,y(t),y^\Delta(t))&=L_{y^\sigma}(t,y^\sigma(t),y^\Delta(t)) \, , \nonumber\\
F_{y^\Delta}(t,y(t),y^\Delta(t))&=L_{y^\sigma}(t,y^\sigma(t),y^\Delta(t))\mu(t)
+L_{y^\Delta}(t,y^\sigma(t),y^\Delta(t))\, ,\nonumber
\end{align}
and the result follows.
\end{proof}


\subsection{The Lagrange problem on time scales}
\label{subsec:Lag}

Now we consider a more general variational problem with
delta-differential side conditions:
\begin{equation}
\label{eq:PrbCO}
\begin{gathered}
J[y(\cdot),u(\cdot)]=\int_{a}^{b}L(t,y(t),u(t))\Delta t \longrightarrow \min \, , \\
y^\Delta(t)=\varphi(t,y(t),u(t)) \, , \\
\left(y(a)=y_a\right), \quad \left(y(b)=y_b\right) \, ,
\end{gathered}
\end{equation}
where $y(\cdot) \in C_{\textrm{rd}}^1[\mathbb{T}]$, $u(\cdot) \in
C_{\textrm{rd}}[\mathbb{T}]$, $y(t)\in\mathbb{R}^n$ and
$u(t)\in\mathbb{R}^m$ for all $t\in \mathbb{T}$, and $m \le n$. We
assume $L: \mathbb{T} \times \mathbb{R}^n \times \mathbb{R}^m
\rightarrow \mathbb{R}$ and $\varphi: \mathbb{T} \times
\mathbb{R}^n \times \mathbb{R}^m \rightarrow \mathbb{R}^n$ to be
$C^1$-functions of $y$ and $u$ for each $t$; and that for each
control function $u(\cdot) \in
C_{\textrm{rd}}[\mathbb{T};\mathbb{R}^m]$ there exists a
correspondent $y(\cdot) \in
C_{\textrm{rd}}^1[\mathbb{T};\mathbb{R}^n]$ solution of the
$\Delta$-differential equation $y^\Delta(t)=\varphi(t,y(t),u(t))$.
We remark that conditions for existence or uniqueness are
available for O$\Delta$E's from the very beginning of the theory
of time scales (see \cite[Theorem~8]{Hilger97}). Roughly speaking,
forward solutions exist, while existence of backward solutions
needs extra assumptions (\textrm{e.g.} regressivity). In control
theory, however, one usually needs only forward solutions, so we
do not need to impose such extra assumptions
\cite{ZbigEwaPaw:CS06}.

We are interested to find necessary conditions for a pair
$\left(y_\ast,u_\ast\right)$ to be a weak local minimizer of $J$.

\begin{defin}
Take an admissible pair $\left(y_\ast,u_\ast\right)$. We say that
$\left(y_\ast,u_\ast\right)$ is a weak local minimizer for
\eqref{eq:PrbCO} if there exist $\delta > 0$ such that
$J[y_\ast,u_\ast] \leq J[y,u]$ for all admissible pairs
$\left(y,u\right)$ satisfying $\|y-y_\ast\|_{1,\infty} +
\|u-u_\ast\|_{\infty} <\delta$.
\end{defin}

\begin{rem}
Problem \eqref{eq:PrbCO} is very general and includes: (i) problem
\eqref{P1} (this is the particular case where $m = n$ and
$\varphi(t,y,u) = u$), (ii) the problem of the calculus of
variations with higher-order delta derivatives \eqref{problema }
(such problem receive special attention in Section~\ref{subsec:HO}
below), (iii) isoperimetric problems on time scales. Suppose that
the isoperimetric condition
\begin{equation*}
I[y(\cdot),u(\cdot)] = \int_a^b g\left(t,y(t),u(t)\right) \Delta t
= \beta \, ,
\end{equation*}
$\beta$ a given constant, is prescribed. We can introduce a new
state variable $y_{n+1}$ defined by
$$y_{n+1}(t)=\int_a^t g(\xi,y(\xi),u(\xi))\Delta \xi,\ t\in\mathbb{T},$$
with boundary conditions $y_{n+1}(a) = 0$ and $y_{n+1}(b) =
\beta$. Then
\begin{equation*}
y_{n+1}^{\Delta}(t)= g\left(t,y(t),u(t)\right),\ t\in\mathbb{T}^k,
\end{equation*}
and we can always recast an isoperimetric problem as a Lagrange
problem \eqref{eq:PrbCO}.
\end{rem}

To establish necessary optimality conditions for \eqref{eq:PrbCO}
is more complicated than for the basic problem of the calculus of
variations on time scales \eqref{eq:EL:B} or \eqref{eq:EL:BSS},
owing to the possibility of existence of abnormal extremals
(Definition~\ref{def:abn}). The abnormal case never occurs for the
basic problem (Proposition~\ref{rem:CV:CN}).

\begin{thm}[The weak maximum principle on time scales]
\label{thm:PMP} If $\left(y_\ast(\cdot),u_\ast(\cdot)\right)$ is a
weak local minimizer of problem \eqref{eq:PrbCO}, then there
exists a set of multipliers $(\psi_{0_\ast}, \psi_\ast(\cdot)) \ne
0$, where $\psi_{0_\ast}$ is a nonnegative constant and
$\psi_\ast(\cdot) : \mathbb{T} \rightarrow \mathbb{R}^n$ is a
delta differentiable function on $\mathbb{T}^k$, such that
$\left(y_\ast(\cdot),u_\ast(\cdot),\psi_{0_\ast},\psi_\ast(\cdot)\right)$
satisfy
\begin{gather}
y_\ast^\Delta(t) = H_{\psi^\sigma}(t,y_\ast(t),u_\ast(t),
\psi_{0_\ast},\psi_\ast^\sigma(t))\, , \quad
\text{($\Delta$-dynamic equation for $y$)}  \label{3}\\
\psi^\Delta_\ast(t)=- H_{y}(t,y_\ast(t),u_\ast(t),
\psi_{0_\ast},\psi_\ast^\sigma(t))\, ,
\quad \text{($\Delta$-dynamic equation for $\psi$)} \label{1} \\
H_{u}(t,y_\ast(t),u_\ast(t),\psi_{0_\ast},\psi_\ast^\sigma(t))=
0\, , \quad \text{($\Delta$-stationary condition)} \label{2}
\end{gather}
for all $t\in\mathbb{T}^k$, where the Hamiltonian
function $H$ is defined by
\begin{equation}
\label{eq:def:Ham} H(t,y,u,\psi_0,\psi^\sigma)= \psi_0 L(t,y,u)
+\psi^\sigma\cdot\varphi(t,y,u) \, .
\end{equation}
If $y(a)$ is free in \eqref{eq:PrbCO}, then
\begin{equation}
\label{eq:trvCondPL:a} \psi_\ast(a) = 0  \, ;
\end{equation}
if $y(b)$ is free in \eqref{eq:PrbCO}, then
\begin{equation}
\label{eq:trvCondPL:b} \psi_\ast(b) = 0  \, .
\end{equation}
\end{thm}

\begin{rem}
From the definition \eqref{eq:def:Ham} of $H$, it follows
immediately that \eqref{3} holds true for any admissible pair
$\left(y(\cdot),u(\cdot)\right)$ of problem \eqref{eq:PrbCO}.
Indeed, condition \eqref{3} is nothing more than the control
system $y_\ast^\Delta(t)=\varphi(t,y_\ast(t),u_\ast(t))$.
\end{rem}

\begin{rem}
For the time scale $\mathbb{T} = \mathbb{Z}$, \eqref{3}-\eqref{2}
reduce to well-known conditions in discrete time (see
\textrm{e.g.} \cite[Ch.~8]{Sethi}): the $\Delta$-dynamic equation
for $y$ takes the form $y(k+1)-y(k) =
H_{\psi}\left(k,y(k),u(k),\psi_0,\psi(k+1)\right)$; the
$\Delta$-dynamic equation for $\psi$ gives $\psi(k+1)-\psi(k) =
-H_{y}\left(k,y(k),u(k),\psi_0,\psi(k+1)\right)$; and the
$\Delta$-stationary condition reads as
$H_u\left(k,y(k),u(k),\psi_0,\psi(k+1)\right) = 0$; with the
Hamiltonian $H = \psi_0 L(k,y(k),u(k)) + \psi(k+1) \cdot
\varphi(k,y(k),u(k))$. For $\mathbb{T} = \mathbb{R}$,
Theorem~\ref{thm:PMP} is known in the literature as \emph{Hestenes
necessary condition}, which is a particular case of the Pontryagin
Maximum Principle \cite{pmp}.
\end{rem}

\begin{cor}[Lagrange multiplier rule on time scales]
\label{cor:LagMultRule} If
$\left(y_\ast(\cdot),u_\ast(\cdot)\right)$ is a weak local
minimizer of problem \eqref{eq:PrbCO}, then there exists a
collection of multipliers $(\psi_{0_\ast}, \psi_\ast(\cdot))$,
$\psi_{0_\ast}$ a~nonnegative constant and $\psi_\ast(\cdot) :
\mathbb{T} \rightarrow \mathbb{R}^n$ a delta differentiable
function on $\mathbb{T}^k$, not all vanishing, such that
$\left(y_\ast(\cdot),u_\ast(\cdot),\psi_{0_\ast},\psi_\ast(\cdot)\right)$
satisfy the Euler-Lagrange equation of the augmented functional
$J^\ast$:
\begin{equation}
\label{eq:prb:Lstar}
\begin{split}
J^\ast[y(\cdot),&u(\cdot),\psi(\cdot)] = \int_{a}^{b}
L^\ast\left(t,y(t),u(t),\psi^\sigma(t),y^\Delta(t)\right) \Delta t \\
&= \int_{a}^{b} \left[ \psi_0 L(t,y(t),u(t)) + \psi^\sigma(t)
\cdot
\left( \varphi(t,y(t),u(t)) - y^\Delta(t) \right)\right] \Delta t \\
&=
\int_{a}^{b}[H(t,y(t),u(t),\psi_0,\psi^\sigma(t))-\psi^\sigma(t)
\cdot y^\Delta (t)]\Delta t \, .
\end{split}
\end{equation}
\end{cor}

\begin{proof}
The Euler-Lagrange equations \eqref{minha:E-L} and \eqref{eq:EL:Boh}
applied to \eqref{eq:prb:Lstar} give
\begin{gather*}
\left(L^\ast_{y^\Delta} -\mu(t)L^\ast_{y}\right)^\Delta =L^\ast_{y} \, , \label{eq:1poss}\\
\left(-\mu(t)L^\ast_{u}\right)^\Delta =L^\ast_{u} \, , \quad
L^\ast_{\psi^\sigma} = 0 \, , \nonumber
\end{gather*}
that is,
\begin{gather}
\left(\psi^\sigma(t) + \mu(t) \cdot H_y\right)^\Delta = - H_{y} \, , \label{eq:1}\\
(-\mu(t)H_{u})^\Delta =H_{u} \, , \label{sei la} \\
y^\Delta(t) = H_{\psi^\sigma}\nonumber \, ,
\end{gather}
where the partial derivatives of $H$ are evaluated at
$(t,y(t),u(t),\psi_0,\psi^\sigma(t))$. Obviously, from \eqref{2}
we obtain \eqref{sei la}. It remains to prove that \eqref{1}
implies \eqref{eq:1} along
$\left(y_\ast(\cdot),u_\ast(\cdot),\psi_{0_\ast},\psi_\ast(\cdot)\right)$.
Indeed, from \eqref{1} we can write $\mu(t) \psi^\Delta(t) = -
\mu(t) H_y$, which is equivalent to $\psi(t) = \psi^\sigma(t) +
\mu(t) H_y$.
\end{proof}

\begin{rem}
Condition \eqref{1} or \eqref{eq:1} imply that along the minimizer
\begin{equation}
\label{eq:new10}
\psi^\sigma(t)=-\int_a^{\sigma(t)}H_{y}(\xi,y(\xi),
u(\xi),\psi_0,\psi^\sigma(\xi))\Delta\xi - c
\end{equation}
for some $c \in \mathbb{R}^n$.
\end{rem}

\begin{rem}
The assertion in Theorem~\ref{thm:PMP} that the multipliers cannot
be all zero is crucial. Indeed, without this requirement, for any
admissible pair $\left(y(\cdot),u(\cdot)\right)$ of
\eqref{eq:PrbCO} there would always exist a set of multipliers
satisfying \eqref{1}-\eqref{2} (namely, $\psi_0 = 0$ and $\psi(t)
\equiv 0$).
\end{rem}

\begin{rem}
Along all the work we consider $\psi$ as a row-vector.
\end{rem}

\begin{rem}
\label{rem:psi0:zoo} If the multipliers
$\left(\psi_0,\psi(\cdot)\right)$ satisfy the conditions of
Theorem~\ref{thm:PMP}, then $\left(\gamma \psi_0,\gamma
\psi(\cdot)\right)$ also do, for any $\gamma
> 0$. This simple observation allow us to conclude that it is
enough to consider two cases: $\psi_0 = 0$ or $\psi_0 = 1$.
\end{rem}

\begin{defin}
\label{def:abn} An admissible quadruple
$\left(y(\cdot),u(\cdot),\psi_0,\psi(\cdot)\right)$ satisfying
conditions \eqref{3}-\eqref{2} (also \eqref{eq:trvCondPL:a} or
\eqref{eq:trvCondPL:b} if $y(a)$ or $y(b)$ are, respectively,
free) is called an extremal for problem \eqref{eq:PrbCO}. An
extremal is said to be normal if $\psi_0 = 1$ and abnormal if
$\psi_0 = 0$.
\end{defin}

So, Theorem~\ref{thm:PMP} asserts that every minimizer is an
extremal.

\begin{prop}
The Lagrange problem on time scales \eqref{eq:PrbCO} has no
abnormal extremals (in particular, all the minimizers are normal)
when at least one of the boundary conditions $y(a)$ or $y(b)$ is
absent (when $y(a)$ or $y(b)$ is free).
\end{prop}

\begin{proof}
Without loss of generality, let us consider $y(b)$ free. We want
to prove that the nonnegative constant $\psi_0$ is nonzero. The
fact that $\psi_0 \ne 0$ follows from Theorem~\ref{thm:PMP}.
Indeed, the multipliers $\psi_0$ and $\psi(t)$ cannot vanish
simultaneously at any point of $t \in \mathbb{T}$. As far as
$y(b)$ is free, the solution to the problem must satisfy the
condition $\psi(b)=0$. The condition $\psi(b)=0$ requires a
nonzero value for $\psi_0$ at $t = b$. But since $\psi_0$ is a
nonnegative constant, we conclude that $\psi_0$ is positive, and
we can normalize it (Remark~\ref{rem:psi0:zoo}) to unity.
\end{proof}

\begin{rem}
\label{rem:CAE} In the general situation abnormal extremals may
occur. More precisely (see proof of Theorem~\ref{thm:PMP}),
abnormality is characterized by the existence of a nontrivial
solution $\psi(t)$ for the system $\psi^\Delta(t) + \psi^\sigma(t)
\cdot \varphi_y = 0$.
\end{rem}

\begin{prop}
\label{rem:CV:CN} There are no abnormal extremals for problem
\eqref{P1}, even in the case $y(a)$ and $y(b)$ are both fixed
($y(a) = y_a$, $y(b) = y_b$).
\end{prop}

\begin{proof}
Problem \eqref{P1} is the particular case of \eqref{eq:PrbCO} with
$y^\Delta(t)= u(t)$. If $\psi_0 = 0$, then the Hamiltonian
\eqref{eq:def:Ham} takes the form $H = \psi^\sigma \cdot\ u$. From
Theorem~\ref{thm:PMP}, $\psi^\Delta = 0$ and $\psi^\sigma = 0$,
for all $t \in \mathbb{T}^k$. Since $\psi^\sigma = \psi + \mu(t)
\psi^\Delta$, this means that $\psi_0$ and $\psi$ would be both
zero, which is not a possibility.
\end{proof}

\begin{cor}
For problem \eqref{P1}, Theorem~\ref{thm:PMP} gives
Theorem~\ref{thm:1}.
\end{cor}

\begin{proof}
For problem \eqref{P1} we have $\varphi(t,y,u) = u$. From
Proposition~\ref{rem:CV:CN}, the Hamiltonian becomes
$H(t,y,u,\psi_0,\psi^\sigma)=L(t,y,u) +\psi^\sigma\cdot u$. By the
$\Delta$-stationary condition \eqref{2} we may write $L_u(t,y,u)
+\psi^\sigma=0$. Now apply \eqref{eq:new10} and the result
follows.
\end{proof}

To prove Theorem~\ref{thm:PMP} we need the following result:

\begin{lem}[Fundamental lemma of the calculus of variations on time scales]
\label{lema fundamental} Let $g\in C_{\textrm{rd}}$,
$g:\mathbb{T}^k\rightarrow\mathbb{R}^n$. Then,
$$
\int_{a}^{b}g(t) \cdot \eta(t)\Delta t=0  \quad \mbox{for all }
\eta\in C_{rd}
$$
if and only if
$$
g(t)=0 \quad \mbox{on }\ \mathbb{T}^k \, .
$$
\end{lem}
\begin{proof}
If $g(t)=0$ on $\mathbb{T}^k$, then obviously $\int_{a}^{b}g(t)
\cdot \eta(t)\Delta t=0$, for all $\eta\in C_{rd}$.

Now, suppose (without loss of generality) that $g(t_0)>0$ for some
$t_0\in\mathbb{T}^k$. We will divide the proof in two steps:

Step 1: Assume that $t_0$ is right scattered. Define in
$\mathbb{T}^k$
\[ \eta(t) = \left\{ \begin{array}{ll}
1 & \mbox{if $t = t_0$};\\
0 & \mbox{if $t \neq t_0$}.\end{array} \right. \] Then $\eta$ is
rd-continuous and
$$\int_a^b g(t)\eta(t)\Delta t=\int_{t_0}^{\sigma(t_0)} g(t)\eta(t)\Delta t=\mu(t_0)g(t_0)>0,$$
which is a contradiction.

Step 2: Suppose that $t_0$ is right dense. Since $g$ is
rd-continuous, then it is continuous at $t_0$. So there exist
$\delta>0$ such that for all
$t\in(t_0-\delta,t_0+\delta)\cap\mathbb{T}^k$ we have $g(t)>0$.

If $t_0$ is left-dense, define in $\mathbb{T}^k$
\[ \eta(t) = \left\{ \begin{array}{ll}
(t-t_0+\delta)^2(t-t_0-\delta)^2 & \mbox{if $t \in (t_0-\delta,t_0+\delta)$};\\
0 & \mbox{otherwise}.\end{array} \right. \] It follows that $\eta$
is rd-continuous and
$$\int_a^b g(t)\eta(t)\Delta t=\int_a^{t_0-\delta} g(t)\eta(t)\Delta t
+\int_{t_0-\delta}^{t_0+\delta} g(t)\eta(t)\Delta t
+\int_{t_0+\delta}^{b} g(t)\eta(t)\Delta t>0,$$
which is a contradiction.

If $t_0$ is left-scattered, define in $\mathbb{T}^k$
\[ \eta(t) = \left\{ \begin{array}{ll}
(t-t_0-\delta)^2 & \mbox{if $t \in [t_0,t_0+\tilde{\delta})$};\\
0 & \mbox{otherwise},\end{array} \right. \] where
$0<\tilde{\delta}<\min\{\mu(\rho(t_0),\delta)\}$. We have: $\eta$
is rd-continuous and
$$\int_a^b g(t)\eta(t)\Delta t=\int_{t_0}^{t_0+\tilde{\delta}} g(t)\eta(t)\Delta t>0,$$
that again leads us to a contradiction.
\end{proof}

\begin{proof} (of Theorem~\ref{thm:PMP})
We begin by noting that $u(t) = \left(u_1(t),\ldots,u_m(t)\right)$
in problem \eqref{eq:PrbCO}, $t\in\mathbb{T}^k$, are arbitrarily
specified functions (controls). Once fixed $u(\cdot) \in
C_{\textrm{rd}}[\mathbb{T}; \mathbb{R}^m]$, then $y(t) =
\left(y_1(t),\ldots,y_n(t)\right)$, $t\in\mathbb{T}^k$, is
determined from the system of delta-differential equations
$y^\Delta(t)=\varphi(t,y(t),u(t))$ (and boundary conditions, if
present). As far as $u(\cdot)$ is an arbitrary function,
variations $\omega(\cdot) \in C_{\textrm{rd}}[\mathbb{T};
\mathbb{R}^m]$ for $u(\cdot)$ can also be considered arbitrary.
This is not true, however, for the variations $\eta(\cdot) \in
$\textrm{C}$_{\textrm{rd}}^1[\mathbb{T}; \mathbb{R}^n]$ of
$y(\cdot)$. Suppose that $(y_\ast(\cdot),u_\ast(\cdot))$ is a weak
local minimizer of $J[\cdot,\cdot]$. Let $\varepsilon \in
(-\delta,\delta)$ be a small real parameter and $y_\varepsilon(t)
= y_\ast(t) + \varepsilon \eta(t)$ (with $\eta(a) = 0$ if
$y(a)=y_a$ is given; $\eta(b) = 0$ if $y(b)=y_b$ is given) be the
trajectory generated by the control $u_\varepsilon(t) = u_\ast(t)
+ \varepsilon \omega(t)$, $\omega(\cdot) \in
C_{\textrm{rd}}[\mathbb{T}; \mathbb{R}^m]$:
\begin{equation}
\label{eq:CSeps}
y_\varepsilon^\Delta(t)=\varphi(t,y_\varepsilon(t),u_\varepsilon(t))
\, ,
\end{equation}
$t \in \mathbb{T}^k$, $\left(y_\varepsilon(a) = y_a\right)$,
$\left(y_\varepsilon(b) = y_b\right)$.  We define the following
function:
\begin{equation*}
\begin{split}
\Phi(\varepsilon) &= J\left[y_\varepsilon(\cdot),
u_\varepsilon(\cdot)\right] = J\left[y_\ast(\cdot)+\varepsilon
\eta(\cdot),u_\ast(\cdot) + \varepsilon \omega(\cdot) \right] \\
&= \int_a^b L\left(t,y_\ast(t)+\varepsilon \eta(t), u_\ast(t) +
\varepsilon \omega(t) \right) \Delta t \, .
\end{split}
\end{equation*}
It follows that $\Phi : (-\delta,\delta) \rightarrow \mathbb{R}$
has a minimum for $\varepsilon=0$, so we must have $\Phi'(0) = 0$.
From this condition we can write that
\begin{equation}
\label{eq:VL} \int_a^b \left[ \psi_0
L_y\left(t,y_\ast(t),u_\ast(t)\right) \cdot \eta(t) + \psi_0
L_{u}\left(t,y_\ast(t),u_\ast(t)\right) \cdot \omega(t) \right]
\Delta t = 0
\end{equation}
for any real constant $\psi_0$. Differentiating \eqref{eq:CSeps}
with respect to $\varepsilon$, we get
\begin{equation*}
\eta^\Delta(t) = \varphi_y(t,y_\varepsilon(t),u_\varepsilon(t))
\cdot \eta(t) + \varphi_{u}(t,y_\varepsilon(t),u_\varepsilon(t))
\cdot \omega(t) \, .
\end{equation*}
In particular, with $\varepsilon = 0$,
\begin{equation}
\label{eq:VSC} \eta^\Delta(t) = \varphi_y(t,y_\ast(t),u_\ast(t))
\cdot \eta(t) + \varphi_{u}(t,y_\ast(t),u_\ast(t)) \cdot \omega(t)
\, .
\end{equation}
Let $\psi(\cdot) \in $\textrm{C}$_{\textrm{rd}}^1[\mathbb{T};
\mathbb{R}^n]$ be (yet) an unspecified function. Multiplying
\eqref{eq:VSC} by $\psi^\sigma(t) =
\left[\psi_1^\sigma(t),\ldots,\psi_n^\sigma(t)\right]$, and
delta-integrating the result with respect to $t$ from $a$ to $b$,
we get that
\begin{equation}
\label{eq:MulPsi:Int} \int_a^b \psi^\sigma(t) \cdot \eta^\Delta(t)
\Delta t = \int_a^b \left[\psi^\sigma(t) \cdot \varphi_y \cdot
\eta(t) + \psi^\sigma(t) \cdot \varphi_{u} \cdot \omega(t)\right]
\Delta t
\end{equation}
for any $\psi(\cdot) \in \textrm{C}_{\textrm{rd}}^1[\mathbb{T};
\mathbb{R}^n]$. Integrating by parts
(see Lemma~\ref{integracao:partes}, formula 1),
\begin{equation}
\label{eq:fzIP}
\begin{split}
\int_a^b \psi^\sigma(t) \cdot \eta^\Delta(t) \Delta t &=
\left.\psi(t) \cdot \eta(t)\right|_a^b
- \int_a^b \psi^\Delta(t) \cdot \eta(t) \Delta t \, ,
\end{split}
\end{equation}
and we can write from \eqref{eq:VL}, \eqref{eq:MulPsi:Int} and
\eqref{eq:fzIP} that
\begin{multline}
\label{eq:LogV} \int_a^b \Bigl[ \left(\psi^\Delta(t) + \psi_0 L_y
+ \psi^\sigma(t) \cdot \varphi_y\right) \cdot \eta(t) \\
+ \left(\psi_0 L_{u} + \psi^\sigma(t) \cdot \varphi_{u}\right)
\cdot \omega(t) \Bigr] \Delta t - \left.\psi(t) \cdot
\eta(t)\right|_a^b = 0
\end{multline}
hold for any $\psi(t)$. Using the definition \eqref{eq:def:Ham} of $H$,
we can rewrite \eqref{eq:LogV} as
\begin{equation}
\label{eq:QF} \int_a^b \left[ \left(\psi^\Delta(t) + H_y\right)
\cdot \eta(t) + H_{u} \cdot \omega(t) \right] \Delta t -
\left.\psi(t) \cdot \eta(t)\right|_a^b = 0 \, .
\end{equation}
It is, however, not possible to employ (yet) Lemma~\ref{lema
fundamental} due to the fact that the variations $\eta(t)$ are not
arbitrary. Now choose $\psi(t) = \psi_\ast(t)$ so that the
coefficient of $\eta(t)$ in \eqref{eq:QF} vanishes:
$\psi_\ast^\Delta(t) = - H_y$ (and $\psi_\ast(a) = 0$ if $y(a)$ is
free, \textrm{i.e.} $\eta(a) \ne 0$; $\psi_\ast(b) = 0$ if $y(b)$
is free, \textrm{i.e.} $\eta(b) \ne 0$). In the normal case
$\psi_\ast(t)$ is determined by
$\left(y_\ast(\cdot),u_\ast(\cdot)\right)$, and we choose
$\psi_{0_\ast} = 1$. The abnormal case is characterized by the
existence of a non-trivial solution $\psi_\ast(t)$ for the system
$\psi_\ast^\Delta(t) + \psi_\ast^\sigma(t) \cdot \varphi_y = 0$:
in that case we choose $\psi_{0_\ast} = 0$ in order to the first
coefficient of $\eta(t)$ in \eqref{eq:LogV} or \eqref{eq:QF} to
vanish. Given this choice of the multipliers, the necessary
optimality condition \eqref{eq:QF} takes the form
\begin{equation*}
\int_a^b H_{u} \cdot \omega(t)  \Delta t = 0 \, .
\end{equation*}
Since $\omega(t)$ can be arbitrarily assigned for all
$t\in\mathbb{T}^k$, it follows from Lemma~\ref{lema fundamental}
that $H_{u} = 0$.
\end{proof}


\subsection{The higher-order problem on time scales}
\label{subsec:HO}

As a corollary of Theorem~\ref{thm:PMP} we obtain the
Euler-Lagrange equation for problem \eqref{problema }.
We first introduce some notation:
\begin{align}
y^0(t)&=y(t),\nonumber\\
y^1(t)&=y^\Delta(t),\nonumber\\
& \ \ \vdots\nonumber\\
y^{r-1}(t)&=y^{\Delta^{r-1}}(t),\nonumber\\
u(t)&=y^{\Delta^r}(t).\nonumber
\end{align}

\begin{thm}
\label{thm:HO:E-L:TS} If $y_\ast\in\mathrm{C}_{rd}^r[\mathbb{T}]$
is a weak local minimizer for the higher-order problem
\eqref{problema }, then
\begin{equation}\label{111}
\psi_\ast^{r-1}(\sigma(t))= - L_{u}(t,x_\ast(t),u_\ast(t))
\end{equation}
holds for all $t\in\mathbb{T}^{k^r}$, where $x_\ast(t) =
\left(y_\ast(t),y_\ast^\Delta(t),\ldots,y_\ast^{\Delta^{r-1}}(t)\right)$
and $\psi_\ast^{r-1}(\sigma(t))$ is defined recursively by
\begin{align}
\psi_\ast^0(\sigma(t))&=-\int_a^{\sigma(t)}L_{y^0}(\xi,x_\ast(\xi),u_\ast(\xi))\Delta\xi+c_0 \, ,\label{222}\\
\psi_\ast^i(\sigma(t))&=-\int_a^{\sigma(t)}\left[L_{y^i}(\xi,x_\ast(\xi),u_\ast(\xi))
+\psi_\ast^{i-1}(\sigma(\xi))\right]\Delta\xi+c_i\label{333},\
i=1,\ldots,r-1 \, ,
\end{align}
with $c_j$, $j = 0,\ldots, r- 1$, constants. If
$y^{\Delta^i}(\alpha)$ is free in \eqref{problema } for some $i
\in \{0,\ldots,r-1\}$, $\alpha \in \{a,\rho^{r-1}(b)\}$, then the
correspondent condition $\psi_\ast^i(\alpha) = 0$ holds.
\end{thm}

\begin{rem}
From \eqref{111}, \eqref{222} and \eqref{333} it follows that
\begin{equation}\label{final}
L_{u}+\sum_{i=0}^{r-1}(-1)^{r-i}\int_a^{\sigma}\cdots\int_a^\sigma
L_{y^i}+[c_i]_{r-i-1}=0,
\end{equation}
where $[c_i]_{r-i-1}$ means that the constant is free from the
composition of the $r-i$ integrals when $i=r-1$ (for simplicity,
we have omitted the arguments in $L_{u}$ and $L_{y^i}$).
\end{rem}

\begin{rem}
If we delta differentiate \eqref{final} $r$ times, we obtain the
delta differentiated equation for the problem of the calculus of
variations with higher order delta derivatives. However, as
observed in Remark~\ref{rem1}, one can only expand formula
\eqref{final} under suitable conditions of delta differentiability
of $\mu(t)$.
\end{rem}

\begin{rem}
For the particular case with $\varphi(t,y,u) = u$, equation
\eqref{eulerint} is \eqref{final} with $r=1$.
\end{rem}

\begin{prop}
\label{prop:NoAbnCaseHO} The higher-order problem on time scales
\eqref{problema } does not admit abnormal extremals, even when the
boundary conditions $y^{\Delta^i}(a)$ and
$y^{\Delta^i}(\rho^{r-1}(b))$, $i = 0,\ldots,r-1$, are all fixed.
\end{prop}

\begin{rem}
\label{rem:Pneeds2rp1points} We require the time scale
$\mathbb{T}$ to have at least $2r+1$ points. Let us consider
problem \eqref{problema } with all the boundary conditions fixed.
Due to the fact that we have $r$ delta derivatives, the boundary
conditions $y^{\Delta^i}(a)=y_a^{i}$ and
$y^{\Delta^i}(\rho^{r-1}(b))=y_b^{i}$ for all $i\in\{0, \ldots,
r-1\}$, imply that we must have at least $2r$ points in order to
have the problem well defined. If we had only this number of
points, then the time scale could be written as
$\mathbb{T}=\{a,\sigma(a),\ldots,\sigma^{2r-1}(a)\}$ and
\begin{equation}
\label{snormal}
\begin{aligned}
\int_{a}^{\rho^{r-1}(b)}
L(t,y(t),&y^\Delta(t),\ldots,y^{\Delta^r}(t))\Delta t \\
&=\sum_{i=0}^{r-1}\int_{\sigma^i(a)}^{\sigma^{i+1}(a)}L(t,
y(t),y^\Delta(t),\ldots,y^{\Delta^r}(t))\Delta t \\
&=\sum_{i=0}^{r-1}L(\sigma^i(a),y(\sigma^i(a)),
y^\Delta(\sigma^i(a)),\ldots,y^{\Delta^r}(\sigma^i(a))),
\end{aligned}
\end{equation}
where we have used the fact that
$\rho^{r-1}(\sigma^{2r-1}(a))=\sigma^r(a)$. Now, having in mind
the boundary conditions and the formula
$$f^\Delta(t)=\frac{f(\sigma(t))-f(t)}{\mu(t)},$$
we can conclude that the sum in \eqref{snormal} would be constant
for every admissible function $y(\cdot)$ and we wouldn't have
nothing to minimize.
\end{rem}

The following technical result is used in the proof of
Proposition~\ref{prop:NoAbnCaseHO}.

\begin{lem}
\label{lemtecn} Suppose that a function
$f:\mathbb{T}\rightarrow\mathbb{R}$ is such that $f^\sigma(t)=0$
for all $t\in\mathbb{T}^k$. Then, $f(t)=0$ for all
$t\in\mathbb{T}\backslash \{a\}$ if $a$ is right-scattered.
\end{lem}

\begin{proof}
First note that, since $f^\sigma(t)=0$, then $f^\sigma(t)$ is
delta differentiable, hence continuous for all $t\in\mathbb{T}^k$.
Now, if $t$ is right-dense, the result is obvious. Suppose that
$t$ is right-scattered. We will analyze two cases: (i) if $t$ is
left-scattered, then $t\neq a$ and by hypothesis
$0=f^\sigma(\rho(t))=f(t)$; (ii) if $t$ is left-dense, then
$f(t)=\lim_{s\rightarrow t^-}f^\sigma(s)=f^\sigma(t)$, by the
continuity of $f^\sigma$. The proof is done.
\end{proof}

\begin{proof} (of Proposition~\ref{prop:NoAbnCaseHO})
Suppose that $\psi_0=0$. With the notation \eqref{eq:PHO:CO}
introduced below, the higher order problem \eqref{problema } would
have the abnormal Hamiltonian given by
$$
H(t,y^0,\ldots,y^{r-1},u,\psi^0,\ldots,\psi^{r-1})
=\sum_{i=0}^{r-2}\psi^{i}(\sigma(t))\cdot y^{i+1}(t)
+\psi^{r-1}(\sigma(t))\cdot u(t)
$$
(compare with the normal Hamiltonian \eqref{eq:normal:Ham:P}).
From Theorem~\ref{thm:PMP}, we can write the system of equations:
\begin{equation}
\label{eq:syst:ab}
\left\{ \begin{array}{ll}
\hat{\psi}^0(t)&=0 \\
\hat{\psi}^1(t)&=-\psi^0(\sigma(t))\\
&\vdots\\
\hat{\psi}^{r-1}(t)&=-\psi^{r-2}(\sigma(t))\\
\psi^{r-1}(\sigma(t))&=0,
\end{array} \right.
\end{equation}
for all $t\in\mathbb{T}^{k^r}$, where we are using the notation
$\hat{\psi}^i(t)={\psi^i}^\Delta(t)$, $i = 0, \ldots, r-1$. From
the last equation, and in view of Lemma~\ref{lemtecn}, we have
$\psi(t)=0$, $\forall t\in\mathbb{T}^{k^{r+1}}\backslash\{a\}$ if
$a$ is right-scattered. This implies that $\hat{\psi}^{r-1}(t)=0$,
$\forall t\in\mathbb{T}^{k^{r}}\backslash\{a\}$ and consequently
$\psi^{r-2}(\sigma(t))=0$, $\forall
t\in\mathbb{T}^{k^{r}}\backslash\{a\}$. Like we did before,
$\psi^{r-2}(t)=0$, $\forall
t\in\mathbb{T}^{k^{r+1}}\backslash\{a,\sigma(a)\}$ if $\sigma(a)$
is right-scattered. Repeating this procedure, we will finally have
$\hat{\psi}^1(t)=0$, $\forall
t\in\mathbb{T}^{k^{r}}\backslash\{a,\ldots,\sigma^{r-2}(a)\}$ if
$\sigma^{i}(a)$ is right-scattered for all $i\in\{0,\ldots,r-2\}$.
Now, the first and second equations in the system
\eqref{eq:syst:ab} imply that $\forall t\in
A=\mathbb{T}^{k^{r}}\backslash\{a,\ldots,\sigma^{r-2}(a)\}$
$$
0=\hat{\psi}^1(t)=-\psi^0(\sigma(t))
=\psi^0(t)+\mu(t)\psi^\Delta(t)=\psi^0(t)\ .
$$
We pick again the first equation to point out that $\psi^0(t)=c$,
$\forall t\in\mathbb{T}^{k^{r+1}}$ and some constant $c$. Since
the time scale has at least $2r+1$ points
(Remark~\ref{rem:Pneeds2rp1points}), the set $A$ is nonempty and
therefore $\psi^0(t)=0,\ \forall t\in\mathbb{T}^{k^{r+1}}$.
Substituting this in the second equation, we get
$\hat{\psi}^1(t)=0,\ \forall t\in\mathbb{T}^{k^{r}}$. As before,
it follows that $\psi^1(t)=d$, $\forall t\in\mathbb{T}^{k^{r+1}}$
and some constant $d$. But we have seen that there exists some
$t_0$ such that $\psi^1(t_0)=0$, hence $\psi^1(t)=0$, $\forall
t\in\mathbb{T}^{k^{r+1}}$. Repeating this procedure, we conclude
that for all $i\in\{0,\ldots,r-1\}$, $\psi^i(t)=0$ at
$t\in\mathbb{T}^{k^{r}}$. This is in contradiction with
Theorem~\ref{thm:PMP} and we conclude that $\psi_0 \ne 0$.
\end{proof}

\begin{proof} (of Theorem~\ref{thm:HO:E-L:TS})
Denoting $\hat{y}(t)=y^\Delta(t)$, then problem \eqref{problema }
takes the following form:
\begin{equation}
\label{eq:PHO:CO}
\begin{gathered}
\mathcal{L}[y(\cdot)]=\int_{a}^{\rho^{r-1}(b)}L(t,y^0(t),y^1(t),
\ldots,y^{r-1}(t),u(t))\Delta t\longrightarrow\min, \\
 \left\{ \begin{array}{l}
\hat{y}^0=y^1 \\
\hat{y}^1=y^2\\
\ \ \ \ \ \vdots\\
\hat{y}^{r-2}=y^{r-1}\\
\hat{y}^{r-1}=u
\end{array} \right.
\end{gathered}
\end{equation}
$$\left(y^i(a)=y_a^i\right),\ \left(y^i\left(\rho^{r-1}(b)\right)=y_b^i\right),\
i=0,\ldots,r-1,\ y_a^i\ \mbox{and}\ y_b^i\in\mathbb{R}^n.$$ System
\eqref{eq:PHO:CO} can be written in the form $y^\Delta = A y + B
u$, where
\begin{equation*}
y = \left(y^0,y^1,\ldots,y^{r-1}\right)
= \left(y_1^0,\ldots,y_n^0,y_1^1,\ldots,y_n^1,\ldots,y_n^{r-1}\right) \in \mathbb{R}^{n r}
\end{equation*}
and the matrices $A$ ($n r$ by $n r$) and $B$ ($n r$ by $n$) are
\begin{equation*}
A = \left(%
\begin{array}{ccccc}
  0 & I & 0 & \cdots & 0 \\
  0 & 0 & I & \cdots & 0 \\
  \vdots & \vdots & \vdots & \ddots &  \vdots \\
  0 & 0 & 0 & \cdots & I \\
  0 & 0 & 0 & \cdots & 0 \\
\end{array}%
\right) \, , \quad
B = col\{0,\ldots,0,I\}
\end{equation*}
in which $I$ denotes the $n$ by $n$ identity matrix, and $0$ the
$n$ by $n$ zero matrix. From Proposition~\ref{prop:NoAbnCaseHO}
we can fix $\psi_0 = 1$: problem \eqref{eq:PHO:CO}
is a particular case of \eqref{eq:PrbCO}
with the Hamiltonian given by
\begin{multline}
\label{eq:normal:Ham:P}
H(t,y^0,\ldots,y^{r-1},u,\psi^0,\ldots,\psi^{r-1})\\
=L(t,y^0,\ldots,y^{r-1},u)+\sum_{i=0}^{r-2}\psi^i(\sigma(\cdot))\cdot
y^{i+1} +\psi^{r-1}(\sigma(\cdot))\cdot u.
\end{multline}
From \eqref{eq:new10} and \eqref{2}, we obtain
\begin{align}
\psi^i(\sigma(t))&=-\int_a^{\sigma(t)}H_{y^i}(\xi,x(\xi),u(\xi),\psi^{\sigma}(\xi))\Delta\xi
+c_i,\ i\in\{0,\ldots,r-1\}\label{11}\\
0&=H_{u}(t,x(t),u(t),\psi^\sigma(t))\label{22},
\end{align}
respectively. Equation \eqref{22} is equivalent to \eqref{111},
and from \eqref{11} we get \eqref{222}-\eqref{333}.
\end{proof}


\section{An example}
\label{subsec:appl}

We end with an application of our higher-order Euler-Lagrange
equation \eqref{final} to the time scale
$\mathbb{T}=[a,b]\cap\mathbb{Z}$, that leads us to the usual and
well-known discrete-time Euler-Lagrange equation (in delta
differentiated form) -- see \textrm{e.g.} \cite{Logan}. Note that
$\forall t\in\mathbb{T}$ we have $\sigma(t)=t+1$ and
$\mu(t)=\sigma(t)-t=1$. In particular, we conclude immediately
that $\mu(t)$ is $r$ times delta differentiable. Also for any
function $g$, $g^\Delta$ exists $\forall t\in\mathbb{T}^k$ (see
Theorem 1.16 (ii) of \cite{livro}) and
$g^\Delta(t)=g(t+1)-g(t)=\Delta g$ is the usual \emph{forward
difference operator} (obviously $g^{\Delta^2}$ exists $\forall
t\in\mathbb{T}^{k^2}$ and more generally $g^{\Delta^r}$ exists
$\forall t\in\mathbb{T}^{k^r}$, $r\in\mathbb{N}$).

Now, for any function $f:\mathbb{T}\rightarrow\mathbb{R}$
and for any $j \in \mathbb{N}$ we have
\begin{align}
\label{eq:exGC}
{\underbrace{\left[\int_a^{\sigma(t)} \left(\int_a^\sigma \cdots \int_a^\sigma f \right)
\Delta \tau\right]}_{j-i\text{ integrals}}}^{\Delta^j}
= f^{\Delta^i \sigma^{j-i}} \, , \quad
i \in \{0,\ldots,j-1\} \, ,
\end{align}
where $f^{\Delta^i \sigma^{j-i}}(t)$ stands for $f^{\Delta^i}(\sigma^{j-i}(t))$.
To see this we proceed by induction. For $j = 1$
\begin{align*}
\int_a^{\sigma(t)}f(\xi)\Delta\xi&=\int_a^{t+1}f(\xi)\Delta\xi
=\int_a^{t}f(\xi)\Delta\xi+\int_t^{t+1}f(\xi)\Delta\xi\\
&=\int_a^{t}f(\xi)\Delta\xi+f(t),
\end{align*}
and then $\left[\int_a^{\sigma(t)}f(\xi)\Delta\xi\right]^\Delta
=f(t)+f^\Delta(t) = f^\sigma$. Assuming that \eqref{eq:exGC} is
true for all $j = 1,\ldots,k$, then
\begin{equation*}
\begin{split}
&{\underbrace{\left[\int_a^{\sigma(t)} \left( \int_a^\sigma \cdots \int_a^\sigma f \right)
\Delta \tau\right]}_{k+1-i\text{ integrals}}}^{\Delta^{k+1}} \\
&= \left(
\underbrace{\int_a^{t} \int_a^\sigma \cdots \int_a^\sigma}_{k+1-i} f
\Delta \tau
+
\underbrace{\int_a^{\sigma(t)} \cdots \int_a^\sigma}_{k-i} f
\Delta \tau
\right)^{\Delta^{k+1}} \\
&=
\left(\underbrace{\int_a^{\sigma(t)} \cdots \int_a^\sigma}_{k-i} f
\Delta \tau
\right)^{\Delta^{k}}
+
\left[\left(\underbrace{\int_a^{\sigma(t)} \cdots \int_a^\sigma}_{k-i} f
\Delta \tau
\right)^{\Delta^{k}}\right]^{\Delta} \\
&= f^{\Delta^i \sigma^{k-i}} + \left(f^{\Delta^i \sigma^{k-i}}\right)^\Delta \\
&= f^{\Delta^i \sigma^{k+1-i}} \, .
\end{split}
\end{equation*}
Delta differentiating $r$ times both sides of equation
\eqref{final} and in view of \eqref{eq:exGC}, we obtain the
Euler-Lagrange equation in delta differentiated form (remember
that $y^0=y$, $\ldots$, $y^{r-1}=y^{\Delta^{r-1}}$,
$y^{\Delta^r}=u$):
\begin{equation*}
L_{y^{\Delta^r}}^{\Delta^r}(t,y,y^\Delta,\ldots,y^{\Delta^r})
+ \sum_{i=0}^{r-1} (-1)^{r-i}
L_{y^{\Delta^{i}}}^{\Delta^i \sigma^{r-i}}(t,y,y^\Delta,\ldots,y^{\Delta^r})
=0.
\end{equation*}


\section{Conclusion}

We introduce a new perspective to the calculus of variations on
time scales. In all the previous works
\cite{Atici06,CD:Bohner:2004,zeidan} on the subject, it is not
mentioned the motivation for having $y^\sigma$ (or $y^\rho$) in
the formulation of problem \eqref{eq:EL:B}. We claim the
formulation \eqref{eq:EL:BSS} without $\sigma$ (or $\rho$) to be
more natural and convenient. One advantage of the approach we are
promoting is that it becomes clear how to generalize the simplest
functional of the calculus of variations on time scales to
problems with higher-order delta derivatives. We also note that
the Euler-Lagrange equation in $\Delta$-integral form
\eqref{eulerint}, for a Lagrangian $L$ with $y$ instead of
$y^\sigma$, follows close the classical condition. Main results of
the paper include: necessary optimality conditions for the
Lagrange problem of the calculus of variations on time scales,
covering both normal and abnormal minimizers; necessary optimality
conditions for problems with higher-order delta derivatives. Much
remains to be done in the calculus of variations and optimal
control on time scales. We trust that our perspective provides
interesting insights and opens new possibilities for further
investigations.


\section*{Acknowledgments}

This work was partially supported by the Portuguese Foundation for
Science and Technology (FCT), through the Control Theory Group
(cotg) of the Centre for Research on Optimization and Control
(CEOC -- \texttt{http://ceoc.mat.ua.pt}). The authors are grateful
to M.~Bohner and S.~Hilger for useful and stimulating comments,
and for them to have shared their expertise on time scales.



\end{document}